\documentclass[twoside,
a4paper]{amsart}
\providecommand{\myincludegraphics}[2][]{\IfFileExists{#2}{\includegraphics[#1]{#2}}{
    File #2 is not found!\typeout{Warning: File #2 is not found!}}}
\usepackage{graphicx}
\PassOptionsToPackage{pdfauthor={Vladimir V. Kisil},%
    pdftitle={Elliptic, Parabolic and Hyperbolic Function Theory--0},%
    pdfsubject={},%
    backref=page,%
  pdfkeywords={symmetries},
  hypertex
}{hyperref}
\IfFileExists{eulervm.sty}{\usepackage{eulervm}
}{}

\usepackage[printedin,extnum
]{myamsart}
\ppnum{LEEDS--MATH--PURE--2004-17
}{arXiv: \href{http://arXiv.org/abs/math.CV/0410399}{math.CV/0410399}}
{2004}

\newcommand{\CPP}{\texttt{C++}}
\newcommand{\NoWEB}{\texttt{noweb}}

\newcommand{\lvec}[1]{\overrightarrow{#1}}
\newcommand{\sperp}{\dashv}

\makeatletter
\def\p@enumi{\thethm.}
\def\myitem#1{\item[(#1)]%
    \protected@edef\@currentlabel
       {\csname p@enumi\endcsname{#1}}%
}
\makeatother

\providecommand{\GiNaC}{\textsf{GiNaC}}

\newcommand{\DSpace}[2]{\ensuremath{ { \dot{\mathbb{#1}}^{#2}} }}
\newcommand{\TSpace}[2]{\ensuremath{ { \widetilde{\mathbb{#1}}^{#2}} }}
\input{mydef}
\renewcommand{\Cliff}[2][\comment]{{\ensuremath{%
\mathcal{C}\kern-0.12em\ell(#1,#2)}}}

\begin{document}
\title[Elliptic, Parabolic and Hyperbolic Function Theory--0]{Elliptic,
  Parabolic and Hyperbolic\\ Analytic Function Theory--0:\\ 
Geometry of Domains}

\author[Vladimir V. Kisil and Debapriya Biswas]%
{\href{http://maths.leeds.ac.uk/~kisilv/}{Vladimir V. Kisil} \and Debapriya Biswas}
\thanks{The first named author is on  leave from the Odessa University.}

\address{%
School of Mathematics\\
University of Leeds\\
Leeds LS2\,9JT\\
UK
}

\email{\href{mailto:kisilv@maths.leeds.ac.uk}{kisilv@maths.leeds.ac.uk}}

\urladdr{\href{http://maths.leeds.ac.uk/~kisilv/}%
{http://maths.leeds.ac.uk/\~{}kisilv/}}

\begin{abstract}
This paper lays down a foundation for a systematic treatment of three
main (elliptic, parabolic and hyperbolic) types of analytic function
theory based on the representation theory of \(\SL\) group. We
describe here geometries of corresponding domains. The principal
r\^ole is played by Clifford algebras of matching types.
\end{abstract}
\AMSMSC{30G35}{22E46}
\maketitle
\tableofcontents
\epigraph{Most attractive feature of most exotic places is the presence
  of a standard tourist accommodation}{}{}
\section{Introduction}
\label{sec:introduction}
Starting from the early age of mathematics as a science we repeatedly
meet the division of various mathematical objects into three main
classes.  In very different areas (equations, quadratic forms,
metrics, manifolds, operators, etc.) these classes preserved names
obtained by the very first example---the classification of conic
sections: elliptic, parabolic, hyperbolic. We will abbreviate this
separation as \emph{EPH-classification}. The \emph{common origin} of
this fundamental division can be seen from the simple picture of a
coordinate line split by the zero into negative and positive
half-axises:
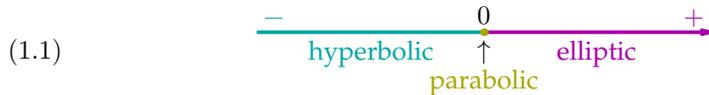
\begin{equation}
  \label{eq:eph-class}
    \raisebox{-15pt}{\myincludegraphics[scale=1]{\jobname.10}}
\end{equation}

However connections between different objects admitting
EPH-classification are not limited to this common source. There are
many deep results linking, for example, ellipticity of quadratic
forms, metrics and operators. On the other hand there are a lot of
white spots and obscure gaps between some subjects as well.

For example,  it is well known that elliptic operators are
effectively treated through complex analysis, which can be naturally
identified as the \emph{elliptic analytic function
  theory}~\cite{Kisil98b,Kisil97c,Kisil01a}. Thus there 
are natural questions about \emph{hyperbolic} and \emph{parabolic}
analytic function theories, which will be of similar importance for
corresponding types of operators. A search for hyperbolic function theory was
initiated in the book~\cite{LavrentShabat77} with some important
advances achieved.

An alternative approach to analytic function theories based on the
representation theory of semisimple Lie groups was developed in the
series of
papers~\cite{Kisil95d,Kisil96d,Kisil98b,Kisil96c,Kisil97c,Kisil94e,Kisil98a,Kisil01a}.
Particularly, a hyperbolic function theory was built
in~\cite{Kisil96d,Kisil98b,Kisil97c} along the same lines as the
elliptic one---standard complex analysis.

This paper makes a further step forwards in this direction. We lay down
foundations for all three (including parabolic!) EPH-types of analytic
function theories. But the present step is rather modest: we  just
study geometries of corresponding domains.
 
\begin{rem}
  Introducing parabolic objects on a common ground with elliptic
  and hyperbolic ones we should warn against two common prejudices
  suggested by 
  picture~\eqref{eq:eph-class}:
\begin{enumerate}
\item \label{it:par-is-zero}
 The parabolic case is unimportant (has ``zero measure'') in
  comparison to the elliptic and hyperbolic ones. As we shall see
  (e.g. Remark~\ref{re:par-more-cayley})  some geometrical
  features are richer in parabolic case.
\item 
  \label{it:par-is-limit} The parabolic case is a limiting situation
  or an intermediate position between the elliptic and hyperbolic: all
  properties of the former can be guessed or obtained as a limit or an
  average from the later two.  Particularly this point of view is
  implicitly supposed in~\cite{LavrentShabat77}.

  Although there are some confirmations of this (e.g.
  Figures~\ref{fig:unit-disks}(\(E\))--(\(H\))), we shall see (e.g.
  Remark~\ref{re:par-is-not-limit}) that some properties of the parabolic
  case cannot be straightforwardly guessed from a combination of
  elliptic and hyperbolic cases.
\end{enumerate}
\end{rem}
An amazing aspect of this topic is a transparent similarity between
all three EPH cases which is combined with some non-trivial exceptions
like \emph{non-invariance} of the upper half plane in the hyperbolic
case (Subsection~\ref{sec:invar-upper-half}) or \emph{non-symmetric}
length and orthogonality in the parabolic case
(Lemma~\ref{item:par-orthogon}). The elliptic case seems to be free
from any such irregularities only because it sets itself the standards to
others.

This paper contains some results and many pictures but almost no
proofs, which are not difficult anyway in most cases. There are only
two notable exceptions (Lemmas~\ref{le:invariance-of-sections}
and~\ref{le:invariance-metric}) when proofs themselves bring
additional insights into the subject.

\section{Elliptic, Parabolic and Hyperbolic Spaces}
\label{sec:ellipt-parab-hyperb}

\subsection[Special linear group and Clifford algebras]{$\SL$ group
  and Clifford Algebras}  
\label{sec:mobi-transf-spr}

We use representations of the \(\SL\) group in Clifford
valued function spaces. There will be three different Clifford algebras
\(\Cliff{e}\), \(\Cliff{p}\), \(\Cliff{h}\)  corresponding to
\emph{elliptic}, \emph{parabolic}, and \emph{hyperbolic} cases
respectively. The notation \(\Cliff{a}\) refers to \emph{any} of these
three algebras.

A Clifford algebra
\(\Cliff{a}\)  as a \(4\)-dimensional  linear space is
spanned by \(1\), \(e_1\), \(e_2\), \(e_1e_2\) with
\emph{non-commutative} multiplication
defined by the identities:
\begin{equation}
  \label{eq:alg-mult}
 e_1^2=-1, \qquad
e_2^2=\left\{
  \begin{array}{cl}
    -1, & \textrm{for \(\Cliff{e}\)---{elliptic} case}\\
    0, & \textrm{for \(\Cliff{p}\)---{parabolic} case}\\
    1, & \textrm{for \(\Cliff{h}\)---{hyperbolic} case}\\
  \end{array}
\right., \qquad 
e_1e_2=-e_2e_1.
\end{equation}
The two-dimensional subalgebra of \(\Cliff{e}\) spanned by \(1\) and
\(\rmi=e_2e_1=-e_1e_2\) is \emph{isomorphic} (and can actually replace
in all calculations!) the field of complex numbers \(\Space{C}{}\).
For any \(\Cliff{a}\) we identify \(\Space{R}{2}\) with the set of
vectors \(w=ue_1+ve_2\), where \((u,v)\in\Space{R}{2}\).  In the elliptic
case of \(\Cliff{e}\) this maps
\begin{equation}
  \label{eq:complexification}
  (u,v)\mapsto e_1(u+\rmi v)=e_1z, \textrm{ with } z=u+\rmi v
  \textrm{ a standard complex number}.
\end{equation} 
We denote \(\Space{R}{2}\) by \(\Space{R}{e}\),
\(\Space{R}{p}\) or \(\Space{R}{h}\) to highlight which of Clifford
algebras is used in the present context. The notation \(\Space{R}{a}\)
assumes \(\Cliff{a}\).

The \(\SL\) group~\cite{HoweTan92,Lang85,MTaylor86} consists of 
\(2\times 2\) matrices
\begin{displaymath}
  \begin{pmatrix}
    a&b\\c&d
  \end{pmatrix},\qquad  \textrm{ with } a, b, c, d\in\Space{R}{}
  \textrm{ and the determinant } ad-bc=1.
\end{displaymath}
An isomorphic realisation of \(\SL\) with the same multiplication is
obtained if we replace a matrix 
\(\begin{pmatrix} 
  a&b\\c&d
\end{pmatrix}
\) by \(\begin{pmatrix}
  a&-b e_1\\c e_1&d
\end{pmatrix}\) within any \(\Cliff{a}\). 
The advantage of the later form is that we can define the \emph{M\"obius 
transformation}  of \(\Space{R}{a}\rightarrow \Space{R}{a}\) for \emph{all}
three algebras \(\Cliff{a}\) by the same expression:
\begin{equation}
  \label{eq:moebius-def}
  \begin{pmatrix}
    a&-b e_1\\c e_1&d
  \end{pmatrix}:\  ue_1+ve_2 \  \mapsto \ 
  \frac {
    a(ue_1+ve_2)-be_1
  } {
    ce_1(ue_1+ve_2)+d
  },
\end{equation}
where the expression \(\frac{a}{b}\) in a non-commutative algebra is
always understood as \(ab^{-1}\), see~\cite{Cnops94a,Cnops02a}. Therefore
\(\frac{ac}{bc}=\frac{a}{b}\) but \(\frac{ca}{cb}\neq\frac{a}{b}\) in
general.

Again in the elliptic case the transformation~\eqref{eq:moebius-def}
is equivalent to 
\begin{displaymath}
  \begin{pmatrix}
    a&-b e_1\\c e_1&d
  \end{pmatrix}:\  e_1z \  \mapsto \ 
  \frac {
    e_1(a(u+e_2e_1 v)-b)
  } {
    -c(u+e_2e_1 v)+d
  }=   e_1\frac {
   az-b
  } {
    -cz+d
  }, \textrm{ where } z=u+\rmi v,
\end{displaymath}
which is the standard form of a M\"obius transformation.
One can straightforwardly verify that the map~\eqref{eq:moebius-def}
is a left action of \(\SL\) on \(\Space{R}{a}\), i.e. \(g_1(g_2
w)=(g_1g_2)w\).

To study the finer structure of M\"obius transformations it is useful to
decompose an element \(g\) of \(\SL\) into the product \(g=g_a g_n g_k\): 
\begin{equation}
  \label{eq:iwasawa-decomp}
  \begin{pmatrix}
    a&-b e_1\\c e_1&d
  \end{pmatrix}=
  {\begin{pmatrix}
      \alpha^{-1} & 0\\0&\alpha
    \end{pmatrix}}
  {\begin{pmatrix}
      1&\chi e_1\\0&1
    \end{pmatrix}}
  {\begin{pmatrix}
      \cos\phi &  e_1\sin\phi\\ 
      e_1\sin\phi & \cos\phi 
    \end{pmatrix}},
\end{equation}
where the values of parameters are as follows:
\begin{equation}
  \label{eq:iwasawa-param}
  \alpha=\sqrt{c^2+d^2}, \qquad
  \chi = \frac{d-a(c^2+d^2)}{c}=\frac{b(c^2+d^2)-c}{d},\qquad
  \phi = \tan^{-1}\frac{c}{d}.
\end{equation} Consequently \(\cos \phi=\frac{d}{\sqrt{c^2+d^2}}\) and
\(\sin \phi=\frac{c}{\sqrt{c^2+d^2}}\).  The
product~\eqref{eq:iwasawa-decomp} gives a realisation of the \emph{Iwasawa
  decomposition }~\cite[\S~III.1]{Lang85} in the form \(\SL=ANK\),
where \(K\) is the maximal compact group, \(N\) is nilpotent and \(A\)
normalises \(N\).

\subsection{Actions of Subgroups and Invariance of Sections} 
\label{sec:mobi-transf}
\begin{figure}[htbp]
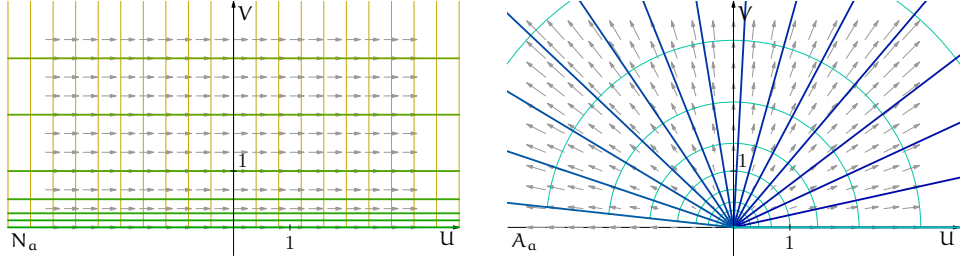

  \myincludegraphics[scale=.75]{\jobname.7}\hfill
   \myincludegraphics[scale=.75]{\jobname.8}
  \caption{Actions of the subgroups \(A\) and \(N\) by M\"obius
    transformations}
  \label{fig:a-n-action}
\end{figure}

In all three EPH cases subgroups the \(A\) and \(N\) act
through M\"obius transformation uniformly: 
\begin{lem} 
\label{le:an-action}
For any type of the Clifford algebra \(\Cliff{a}\):
  \begin{enumerate}
  \item The subgroup \(N\) defines shifts \(ue_1+ve_2 \mapsto
    (u+\chi)e_1+ue_2 \) along the ``real'' axis \(U\) by
    \(\chi\).\\
    The vector field of the derived representation is \(dN_a{(u,v)}=(1,0)\).
  \item The subgroup \(A\) defines dilations \(ue_1+ve_2 \mapsto
    \alpha^2(ue_1+ve_2)\) by the factor \(\alpha^2\) which fixes origin
    \((0,0)\).\\
    The vector field of the derived representation is \(dA_a{(u,v)}=(2u,2v)\).
  \end{enumerate}
\end{lem}
Orbits and vector fields corresponding to the \emph{derived
representation}~\cite[\S~6.3]{Kirillov76}, \cite[Chap.~VI]{Lang85} 
of the Lie algebra \(\algebra{sl}_2\) for subgroups \(A\) and \(N\)
are shown in Figure~\ref{fig:a-n-action}. Thin transverse lines
join points of orbits corresponding to the same values of the parameter. 

By contrast the actions of the subgroup  \(K\) is significantly
different between the EPH cases and  correlates with names chosen for
\(\Cliff{e}\), \(\Cliff{p}\), \(\Cliff{h}\):
\begin{figure}[htbp]
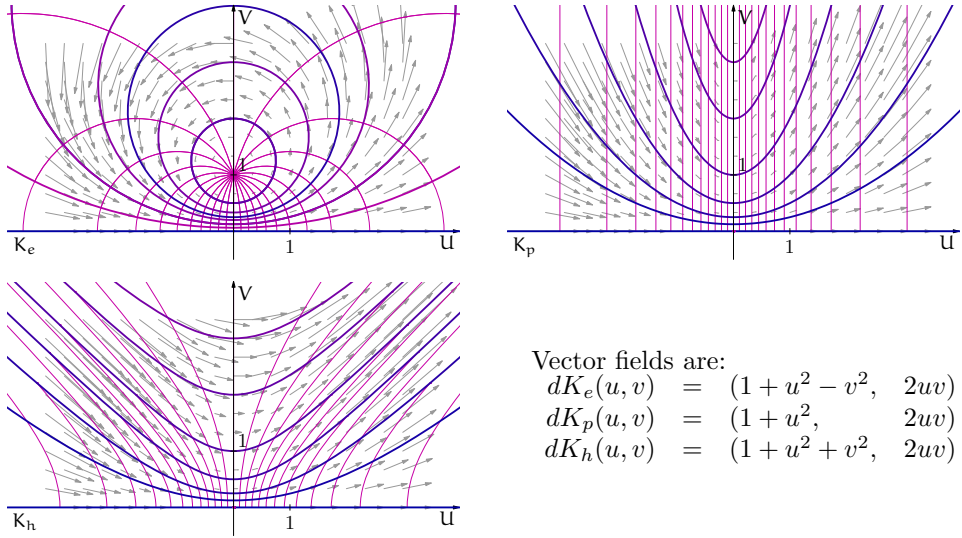

  \centering
\myincludegraphics[scale=.75
]{\jobname.4}\hfill
\myincludegraphics[scale=.75
]{\jobname.5}\\[2mm]
\parbox{.45\textwidth}{
\myincludegraphics[scale=.75
]{\jobname.6}
}
\hfill \parbox{.45\textwidth}{
  Vector fields are:\\
  \(\begin{array}{rcll}
    dK_e(u,v)&=&(1+u^2-v^2,&2uv)\\
    dK_p(u,v)&=&(1+u^2,&2uv)\\
    dK_h(u,v)&=&(1+u^2+v^2,&2uv)\\
  \end{array}\)
}
  \caption{Action of the \(K\) subgroup. The corresponding orbits are
    circles, parabolas and hyperbolas.}
  \label{fig:k-sungroup}
\end{figure}
\begin{lem} 
\label{le:k-action}
The actions of the subgroup \(K\) in three cases are as
  follows: 
  \begin{enumerate}
    \myitem{e} \label{item:circle-desc} For \(\Cliff{e}\) the orbits
    of \(K\) are circles. A circle with centre at \((0,
    (v+v^{-1})/2)\) passing through two points \((0,v)\)
    and \((0,v^{-1})\).\\
    The vector field of the derived representation is
    \(dK_e{(u,v)}=(u^2-v^2+1, 2uv)\).  
    \myitem{p}
    \label{item:parab-desc} For \(\Cliff{p}\) the orbits of \(K\) are
    parabolas with the vertical axis \(V\). A parabola passing through
    \((0,v/2)\) has its horizontal directrix passing through
    \((0,(v-v^{-1}/2)\)
    and focus at \((0,(v+v^{-1})/2)\). \\
    The vector field of the derived representation is
    \(dK_p{(u,v)}=(u^2+1, 2uv)\).  
    \myitem{h} \label{item:hyperb-desc}
    For \(\Cliff{h}\) the orbits of \(K\) are hyperbolas with asymptotes
    parallel to lines \(u=\pm v\). A hyperbola passing through the point
    \((0,v)\) has the focal distance between foci \(2p\), where
    \(p=\frac{v^2+1}{\sqrt{2}v}\) and the upper focus is located at
    \((0,f)\) with:
    \begin{displaymath}
    f=\left\{\begin{array}{ll}
      p-\sqrt{\frac{p^2}{2}-1}, &\textrm{ for } 0<v<1; \textrm{
        and }\\
      p+\sqrt{\frac{p^2}{2}-1}, & \textrm{ for } v\geq 1.
    \end{array}\right.
  \end{displaymath}
  The vector field of the derived representation is
    \(dK_h{(u,v)}=(u^2+v^2+1, 2uv)\).
  \end{enumerate}
\end{lem} 
Orbits and the corresponding derived actions of the subgroup
\(K\) are shown in Figure~\ref{fig:k-sungroup}.
\begin{rem}  
  \begin{enumerate}
  \item The  values of all three vector fields \(dK_e\), \(dK_p\) and \(dK_h\)
    coincide on the ``real'' \(U\)-axis \(v=0\), i.e. they are three
    different extensions into the domain of the same boundary condition.
  \item \label{it:hyp-orbit}
    The hyperbola passing through the point \((0,1)\) has the shortest
    focal length \(\sqrt{2}\) among all other hyperbolic orbits; two
    hyperbolas passing through \((0,v)\) and \((0,v^{-1})\) have the
    same focal length and are related to each other as explained in
    Remark~\ref{it:hyp-object}. 
\end{enumerate}
\end{rem}
\begin{defn}
  We use the word \emph{cycle}
  to  denote straight lines \textbf{and one of the following}
  \begin{enumerate}
  \myitem{e} Circles in the elliptic case;
  \myitem{p} Parabolas with a vertical axis of symmetry in the parabolic
    case;
  \myitem{h} Equilateral hyperbolas with a vertical axis of symmetry in the
    hyperbolic case. 
  \end{enumerate}
  Moreover the words \emph{parabola} and \emph{hyperbola} in this paper
  always assume only ones of the above described types.

  \emph{Centre} of a cycle is its geometrical centre for a circle or
  a hyperbola and its focus for a parabola. Centres of straight lines
  are at infinity.
\end{defn}

Using the Lemmas~\ref{le:an-action} and~\ref{le:k-action} we can give
an easy proof of invariance for corresponding cycles. 

\begin{lem} \label{le:invariance-of-cycles}
  M\"obius transformations preserve the cycles in the
  upper   half plane, i.e.:
  \begin{enumerate}
  \myitem{e} For \(\Cliff{e}\) M\"obius transformations map circles to circles.
  \myitem{p} For \(\Cliff{p}\) M\"obius transformations map parabolas to parabolas.
  \myitem{h} For \(\Cliff{h}\) M\"obius transformations map hyperbolas to hyperbolas. 
  \end{enumerate}
\end{lem}
\begin{proof}
\begin{figure}[htbp]
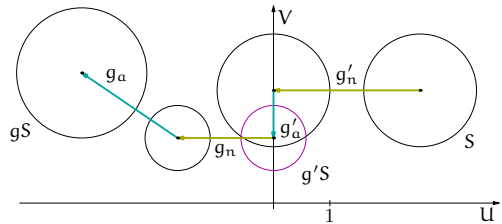

  \centering
  \myincludegraphics[scale=.75]{\jobname.14}
  \caption[Decomposition of an arbitrary Moebius
  transformation]{The decomposition of an arbitrary M\"obius
    transformation \(g\) into a product \(g=g_a g_n g_k g_a' g_n'\).}
  \label{fig:moeb-decomp}
\end{figure}
  Our first observation is that the subgroups \(A\) and \(N\) obviously
  preserve all circles, parabolas, hyperbolas and straight lines in all
  \(\Cliff{a}\). Thus we use subgroups \(A\) and \(N\) to  fit a given
  cycle exactly on a particular orbit of subgroup \(K\) shown on
  Figure~\ref{fig:k-sungroup} of the corresponding type. 

  To this end for an arbitrary cycle \(S\) we can find \(g_n'\in N\)
  which puts centre of \(S\) on the \(V\)-axis, see
  Figure~\ref{fig:moeb-decomp}. Then there is a unique \(g_a'\in A\)
  which scales it exactly to an orbit of \(K\), e.g. for a
  circle passing through point \((0,v_1)\) and \((0,v_2)\) the scaling
  factor is \(\frac{1}{\sqrt{v_1v_2}}\) accordingly to
  Lemma~\ref{item:circle-desc}. Let \(g'=g_a' g_n'\), then for any
  element \(g\in \SL\) using the Iwasawa decomposition of \(g
  g'^{-1}=g_a g_n g_k\) we get the presentation \(g=g_a g_n g_k g_a'
  g_n'\) with \(g_a, g_a'\in A\), \(g_n, g_n'\in N\) and \(g_k \in
  K\).

  Then the image \(g'S\) of the cycle \(S\) under \(g'=g_a' g_n'\)
  is a cycle itself in the obvious way, then \(g_k (g'S)\) is again a
  cycle since \(g'S\) was arranged to coincide with a \(K\)-orbit, and
  finally  \(gS=g_a g_n (g_k (g'S))\) is a cycle due to the obvious
  action of \(g_a g_n\), see Figure~\ref{fig:moeb-decomp} for an
  illustration. 
\end{proof}

\subsection{Lengths and Orthogonality}
\label{sec:lenghts-orth}
The invariance of cycles (see Lemma~\ref{le:invariance-of-cycles})
suggests using them in the r\^ole of circles in each of the EPH cases 
 and play \emph{the standard mathematical game}: turn some
properties of classical objects into definitions of new ones.
\begin{defn}
  \label{de:length}
  The \emph{length} \(l_a(\lvec{AB})\) of a vector \(\lvec{AB}\) in
  \(\Space{R}{a}\) is defined as a real valued function, such that for
  a fixed \(A\) the level curves of \(l_a(\lvec{AB})\) are of
  corresponding shapes: circles, parabolas or hyperbolas.
\end{defn}
\begin{lem} The following are lengths in the sense of Definition~\ref{de:length}:
  \begin{enumerate}
  \myitem{e} In the elliptic case: the Euclidean metric
    \(l_e(ue_1+ve_2)=u^2+v^2\).  
    \myitem{p} In the parabolic case: a monotonic function of focal
    length:
    \begin{equation}
      \label{eq:par-focal-length}
      l_p(ue_1+ve_2)=\sqrt{u^2+v^2}-v
    \end{equation}
    for a parabola with
    focus \(A\) passing through \(B\). 
    Note, that \(l(\lvec{AB})\neq l(\lvec{BA})\)!  
    \myitem{h} In the hyperbolic case: the Minkowski
    metric \(l_h(ue_1+ve_2)={u^2-v^2}\).
\end{enumerate}
\end{lem}
\begin{rem}
  In the elliptic and hyperbolic cases the above lengths are conveniently
  defined by the Clifford algebra multiplication
  \begin{displaymath}
    l_{e,h}(ue_1+ve_2)=-(ue_1+ve_2)^2.    
  \end{displaymath} It will be interesting to find some sort of such
  relation for the parabolic length~\eqref{eq:par-focal-length} as
  well.
\end{rem}
\begin{defn}
  We say that a vector \(\lvec{AB}\) is \emph{s-orthogonal} to a
  vector \(\lvec{CD}\) and denote it \(\lvec{AB} \sperp \lvec{CD}\)
  (for the reasons clear from Lemma~\ref{item:par-orthogon})
  if the function \(l(\lvec{AB}+\epsilon \lvec{CD})\) of a variable
  \(\epsilon\) has a local extremum at \(\epsilon=0\)
  (i.e. orthogonality provides the shortest length).
\end{defn}
Again s-orthogonality turns out to be the usual orthogonality in the
elliptic case. For the two other cases the description is given as
follows:
\begin{lem} 
  A vector \(ue_1+ve_2\) is \emph{s-orthogonal} to a vector
  \(u'e_1+v'e_2\) if in terms  of Euclidean geometry:
  \begin{enumerate} 
  \myitem{e} In the elliptic case: two vectors form a right angle, or
  analytically \(u u'+v v'=0\).
  \myitem{p} \label{item:par-orthogon} In the parabolic case: the vector
  \(ue_1+ve_2\) bisects the angle between \(u'e_1+v'e_2\) and the
  vertical directions or analytically:
  \begin{equation}
    \label{eq:par-s-perp}
    u'u-v'l_p(ue_1+ve_2)=u'u-v'(\sqrt{u^2+v^2}-v)=0.
  \end{equation}
  Note that \( \lvec{AB} \sperp \lvec{CD}\) does
  \textbf{not} necessarily imply \( \lvec{CD} \sperp \lvec{AB} \)!
  \myitem{h}  In the hyperbolic case the angles between two vectors
  are bisected by lines parallel to \(u=\pm v\), or analytically
  \(u' u-v' v=0\). 
  \end{enumerate}
\end{lem}
\begin{rem}
  \label{re:par-is-not-limit}
  If one tries to devise a parabolic length as a limit or an intermediate
  case between the elliptic \(l_e=u^2+v^2\) and hyperbolic \(l_p=u^2-v^2\)
  lengths  then the only possible guess is \(l'_p=u^2\), which  is too
  trivial for an interesting geometry. 
  
  Similarly the only orthogonality conditions bridging elliptic
  \(u_1 u_2+v_1 v_2=0\) and hyperbolic \(u_1 u_2-v_1 v_2=0\) seems to
  be \(u_1 u_2=0\) which is again too trivial. This support our
  Remark~\ref{it:par-is-limit}. 
\end{rem}

\subsection{Zero Radius Cycles, Invariant Measure and Compactification}
\label{sec:zero-radius-cycles}
Of course, M\"obius transformations may not preserve centres
of cycles. However this happens in a trivial way for ``zero radius''
cycles, as follows.
\begin{lem}
  A zero-radius cycle with centre at \((u_0,v_0)\) defined by
  the equation \(l_a(u-u_0,v-v_0)=0\) is:
\begin{enumerate}
\myitem{e} a single point \((u_0,v_0)\) in the elliptic case;
\myitem{p} \label{it:par-zero-rad}
 the vertical upward directed ray with origin at \((u_0,v_0)\)
in the parabolic case; 
\myitem{h} the \emph{light cone} with origin at  \((u_0,v_0)\)
defined by the equation 
\begin{equation}
  \label{eq:light-cone}
   \qquad (u-u_0)^2-(v-v_0)^2=0
\end{equation}
 in the hyperbolic case.
\end{enumerate}
\end{lem}
In the elliptic and hyperbolic cases it is often
useful~\cite{Cnops94a,Cnops02a} to consider zero radius cycles
instead of corresponding points which is known as
Fillmore-Springer-Cnops construction. The same advantages are expected
in the parabolic case as well.  Among many useful applications the
embedding of \(\Space{R}{a}\) into a bigger space of spheres produces
the invariant measure in an elegant way~\cite{Cnops94a,Cnops02a}. We
give another proof based on the Iwasawa decomposition.
\begin{lem}
\label{le:invariance-metric}
  A M\"obius invariant measure on \(\Space{R}{a}\) is given by \(\frac{du\,dv}{v^2}\).
\end{lem}
\begin{proof}
  Let \(f(u,v)\,du\,dv\) be an invariant measure. Then considering
  shifts generated by the subgroup \(N\) (see
  Figure~\ref{fig:a-n-action}) we conclude that \(f(u,v)\) is
  independent of \(u\), thus we denote it by \(f(v)\). The dilations
  generated by the subgroup \(A\) (see Figure~\ref{fig:a-n-action})
  put the restriction \(f(v)=cv^{-2}\) which is obviously compatible
  with any \(K\) action since \(\partial_v\) components of all vector
  fields \(dK_a\) are the same.
\end{proof}

Another important r\^ole of the zero radius cycles is the proper
compactification of \(\Space{R}{a}\). Indeed the initial space
\(\Space{R}{a}\) is not a closed set under a generic M\"obius
transformations. In the elliptic case the problem is solved by the
compactification of \(\Space{R}{e}\) with a point \(\infty\) at
infinity. However in the parabolic and hyperbolic cases the singularity of
the M\"obius transform is not localised in a single point---the denominator
vanish for the whole zero radius cycle.  Thus in each EPH case the
correct compactification is made by a \emph{zero radius cycle at
  infinity}. Of course in the elliptic case this is still a point,
but for the two other cases the result is significantly different.

It is common to identify the compactification \(\DSpace{R}{e}\) of the
space \(\Space{R}{e}\) by a point \(\infty\) with a Riemann sphere.
This model can be visualised by the stereographic
projection~\cite[\S~18.1.4]{BergerII}. The projection from the centre
of a sphere provides a model for the compactification
\(\TSpace{R}{p}\) in the parabolic case. The space \(\TSpace{R}{p}\)
is represented by a sphere where all pairs of opposite points are
identified. The ``half of the equator'' in this model represents the
parabolic zero radius cycle (see Lemma~\ref{it:par-zero-rad}) at
infinity. More informative models are provided by the
Fillmore-Springer-Cnops construction, which represent M\"obius
transformations through orthogonal rotations in the bigger space of
spheres~\cite{Cnops94a,Cnops02a}.

The hyperbolic case produces its own caveats. A compactification of the
hyperbolic space \(\Space{R}{h}\) by a light cone at infinity will
produce a closed M\"obius invariant object. However it will not
be satisfactory for some other reasons explained in the next
Subsection.

\subsection{(Non)-Invariance of The Upper Half Plane}
\label{sec:invar-upper-half}
\begin{figure}[htbp]
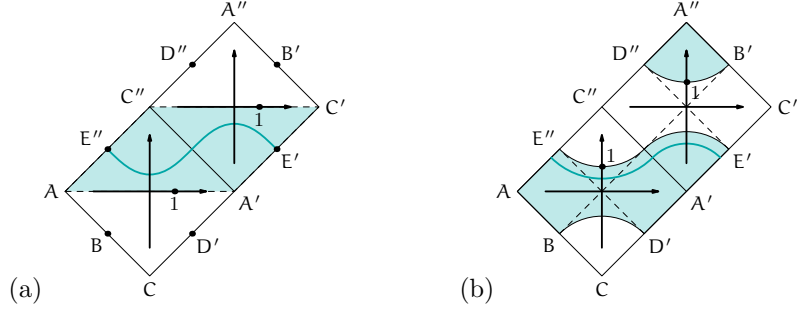

  \centering
     (a)\myincludegraphics[scale=.75]{\jobname.13} \hspace{.1\textwidth}
    (b)\myincludegraphics[scale=.75
]{\jobname.15}
  \caption[Hyperbolic objects in the double cover]{Hyperbolic objects
    in the double cover of \(\Space{R}{h}\):\\ 
  (a) the ``upper'' half plane;\qquad (b) the unit circle.}
  \label{fig:hyp-upper-half-plane}
\end{figure}

The important difference between the hyperbolic case and the two
others is that
\begin{lem}
  In the elliptic and parabolic cases the upper halfplane in
  \(\Space{R}{a}\) is preserved by M\"obius transformations from
  \(\SL\). However in the hyperbolic case any point \((u,v)\) with
  \(v>0\) can be mapped into an arbitrary point \((u',v')\) with \(v'\neq 0\).
\end{lem}
The lack of invariance in the hyperbolic case has many important
consequences in seemingly different areas, for example:
\begin{figure}[htbp]
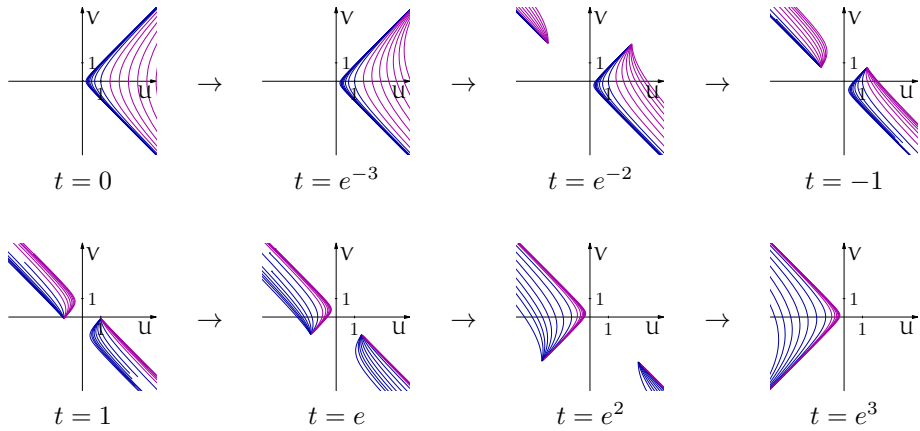

  \parbox[t]{.2\textwidth}{
    \begin{center}
      \myincludegraphics[scale=.70]{\jobname.20}\\
      \(t=0\)
    \end{center}
  } \hfill\raisebox{-1.2cm}{\(\to\)} \hfill
  \parbox[t]{.2\textwidth}{
    \begin{center}
      \myincludegraphics[scale=.70]{\jobname.21}\\
      \(t=e^{-3}\)
    \end{center}
  }   \hfill \raisebox{-1.2cm}{\(\to\)} \hfill
  \parbox[t]{.2\textwidth}{
    \begin{center}
      \myincludegraphics[scale=.70]{\jobname.22}\\
      \(t=e^{-2}\)
    \end{center}
  } \hfill \raisebox{-1.2cm}{\(\to\)} \hfill
    \parbox[t]{.2\textwidth}{
    \begin{center}
      \myincludegraphics[scale=.70]{\jobname.23}\\
      \(t=-1\)
    \end{center}
  }\\[5mm]
  \parbox[t]{.2\textwidth}{
    \begin{center}
      \myincludegraphics[scale=.70]{\jobname.24}\\
      \(t=1\)
    \end{center}
  }  \hfill \raisebox{-1.2cm}{\(\to\)} \hfill
  \parbox[t]{.2\textwidth}{
    \begin{center}
      \myincludegraphics[scale=.70]{\jobname.25}\\
      \(t=e\)
    \end{center}
  } \hfill \raisebox{-1.2cm}{\(\to\)} \hfill
  \parbox[t]{.2\textwidth}{
    \begin{center}
      \myincludegraphics[scale=.70]{\jobname.26}\\
      \(t=e^2\)
    \end{center}
  } \hfill \raisebox{-1.2cm}{\(\to\)} \hfill
  \parbox[t]{.2\textwidth}{
    \begin{center}
      \myincludegraphics[scale=.70]{\jobname.27}\\
      \(t=e^3\)
    \end{center}
  }
  \caption{Eight frames from a continuous transformation from future
    to the past parts   of the light cone.}
  \label{fig:future-to-past}
\end{figure}
\begin{description}
\item[\textbf{Geometry}] \(\Space{R}{h}\) is not split by the real
  axis into two disjoint pieces: there is a continuous path (through the
  light cone at infinity) from the upper half plane to lower which
  does not cross the real axis (see the \(\sin\)-like joined two pieces of
  the hyperbola in Figure~\ref{fig:hyp-upper-half-plane}(a)).
\item[\textbf{Physics}] There is no M\"obius invariant way to separate
  ``past'' and ``future'' parts of the light cone~\cite{Segal76}, i.e.
  there is a continuous family of M\"obius transformations reversing
  the arrow of time. For example, The family of matrices \(
  \begin{pmatrix}
    1&-te_2\\te_2&1
  \end{pmatrix}\), \(t\in [0,\infty)\) provide the transformations and
  Figure~\ref{fig:future-to-past} presents images for eight values of \(t\). 
\item[\textbf{Analysis}] There is no a possibility to split
  \(\FSpace{L}{2}(\Space{R}{})\) space of function into a direct sum
  of the Hardy space of functions having an analytic extension into
  the upper half plane and its non-trivial complement, i.e. any function
  from \(\FSpace{L}{2}(\Space{R}{})\) has an ``analytic extension''
  into the upper half plane,  see~\cite{Kisil97c}. 
\end{description}
 All the above problems can be resolved in the following
way~\cite[\S~A.3]{Kisil97c}.  We take two copies
\(\Space[+]{R}{h}\) and \(\Space[-]{R}{h}\) of \( \Space{R}{h} \),
depicted by squares \(ACA'C''\) and \(A'C'A''C''\) in
Figure~\ref{fig:hyp-upper-half-plane} correspondingly. The boundaries
of these squares are light cones at infinity and we glue
\(\Space[+]{R}{h}\) and \(\Space[-]{R}{h}\) in such a way that the
construction is invariant under the natural action of the M\"obius
transformation.  That is achieved if the same letters \(A\), \(B\), \(C\),
\(D\), \(E\) in Figure~\ref{fig:hyp-upper-half-plane} are
identified regardless the number of attached primes.  This
aggregate denoted by \(\TSpace{R}{h}\) is a two-fold cover of
\(\Space{R}{h}\). The hyperbolic ``upper'' half plane in
\(\TSpace{R}{h}\) consists of the upper halfplane in
\(\Space[+]{R}{h}\) and the lower one in \(\Space[-]{R}{h}\).  A
similar conformally invariant two-fold cover of the Minkowski
space-time was constructed in~\cite[\S~III.4]{Segal76} in connection
with the red shift problem in extragalactic astronomy.

\begin{rem}
  \begin{enumerate}
  \item \label{it:hyp-object}
    The hyperbolic orbit of the \(K\) subgroup in the  \(
    \TSpace{R}{h}\) consists of two branches of the hyperbola passing through
    \((0,v)\) in  \(\Space[+]{R}{h}\) and \((0,-v^{-1})\) in
    \(\Space[-]{R}{h}\), see Figure~\ref{fig:hyp-upper-half-plane}. 
    As explained in Remark~\ref{it:hyp-orbit} they
    both have the same focal length.
  \item  The ``upper'' halfplane is bounded by two disjoint ``real''
    axises denoted by \(AA'\) and \(C'C''\) on Figure~\ref{fig:hyp-upper-half-plane}.
  \end{enumerate}
\end{rem}

For the hyperbolic Cayley transform in the next subsection we 
need the conformal version of the hyperbolic unit disk.
 We define it in \( \TSpace{R}{1,1} \) as follows:
\begin{eqnarray*}
\TSpace{D}{}&=&\{(ue_1+ve_2) \such l_h(ue_1+ve_2)<-1,\ u\in \mathbb{R}^{1,1}_+ \}\\
&&{}\cup \{(ue_1+ve_2) \such l_h(ue_1+ve_2)>-1,\ u\in \mathbb{R}^{1,1}_- \}.
\end{eqnarray*}
It can be shown that \( \TSpace{D}{}\) is conformally invariant and 
has a boundary \(\TSpace{T}{}\)---the two copies of the unit circles in 
\(\Space[+]{R}{1,1}\) and \(\Space[+]{R}{1,1}\). We call \(\TSpace{T}{}\) 
the \emph{(conformal) unit circle} in \(\Space{R}{1,1}\). 
Figure~\ref{fig:hyp-upper-half-plane} illustrates\footnote{Note that
  similar figures in papers~\cite{Kisil97c,Kisil02c} have letters
  \(D'\) and \(E'\) misplaced.} the geometry of the ``upper'' half plane
as well as the conformal unit disk in 
\(\TSpace{R}{1,1}\) conformally equivalent 
to it.

\subsection{The Cayley Transform and Unit ``Circles''}\vspace{-2mm}
\label{sec:unit-circles}

The upper half plane is the universal starting point for an analytic
function theory of any EPH type. However universal models are rarely
best suited to particular circumstances. For many reasons it is more
convenient to consider analytic functions on the unit disk rather than on the
upper half plane, although both theories are completely
isomorphic, of course. This isomorphism is delivered by the \emph{Cayley
  transform}. 

Let \(\sigma=e_2^2\), i.e. \(-1\), \(0\), or \(1\) as
in~\eqref{eq:alg-mult}. Then the first possibility to define the Cayley
transform is given by the matrix \(C=\begin{pmatrix}
  1 &  -e_2 \\ \sigma e_2 & 1
\end{pmatrix}\) with determinant \(1\). It can be applied as the M\"obius
transformation
\begin{equation}
  \label{eq:cayley-points}
  \begin{pmatrix}
  1 & -e_2 \\ \sigma e_2 & 1
\end{pmatrix}:\  w=(ue_1+ve_2) \mapsto Cw=
\frac{(ue_1+ve_2)-e_2}{ \sigma e_2 (ue_1+ve_2)+1}
\end{equation}
to a point \((ue_1+ve_2)\in\Space{R}{a}\).  Alternatively it acts by conjugation
\(g_C=CgC^{-1}\) on an element \(g\in\SL\):
\begin{equation}
  \label{eq:cayley-matr}
  g_C= \frac{1}{2}
  \begin{pmatrix}
    1 & - e_2 \\ \sigma e_2 & 1
  \end{pmatrix}
  \begin{pmatrix}
    a & -be_1 \\ c e_1 & d
  \end{pmatrix}
  \begin{pmatrix}
    1 & e_2 \\ -\sigma e_2 & 1
  \end{pmatrix} 
\end{equation} The connection between the two
forms~\eqref{eq:cayley-points} and~\eqref{eq:cayley-matr} of the Cayley
transform is given by \(g_C Cw= C(gw)\), i.e. \(C\) intertwines
the actions of \(g\) and \(g_C\).

The Cayley transform \((u'e_1+v'e_2)= C(ue_1+ve_2)\) 
in the elliptic case is very
important~\cite[\S~IX.3]{Lang85}, \cite[Ch.~8,
(1.12)]{MTaylor86}. 
The transformation \(g\mapsto g_C\)~\eqref{eq:cayley-matr}
is an isomorphism of the groups \(\SL\) and \(\mathrm{SU}(1,1)\),
namely in \(\Cliff{e}\) we have
\begin{equation}
  \label{eq:cayley-elliptic}
  g_C= 
  \begin{pmatrix}
    f & h \\ -h & f
\end{pmatrix}, \textrm{ with }
f= (a+d)-(c-b)e_2e_1 \textrm{ and } h=(a-d)e_2-(b+c)e_1.
\end{equation}
Under the map \(\Space{R}{e}\to \Space{C}{}\)~\eqref{eq:complexification}
this matrix becomes \(
\begin{pmatrix}
  \alpha & \beta \\ \bar{\alpha} & \bar{\beta}
\end{pmatrix}
\), i.e. the standard form of elements of
\(\mathrm{SU}(1,1)\)~\cite[\S~IX.1]{Lang85}, \cite[Ch.~8,
(1.11)]{MTaylor86}.  

\begin{figure}[ht]
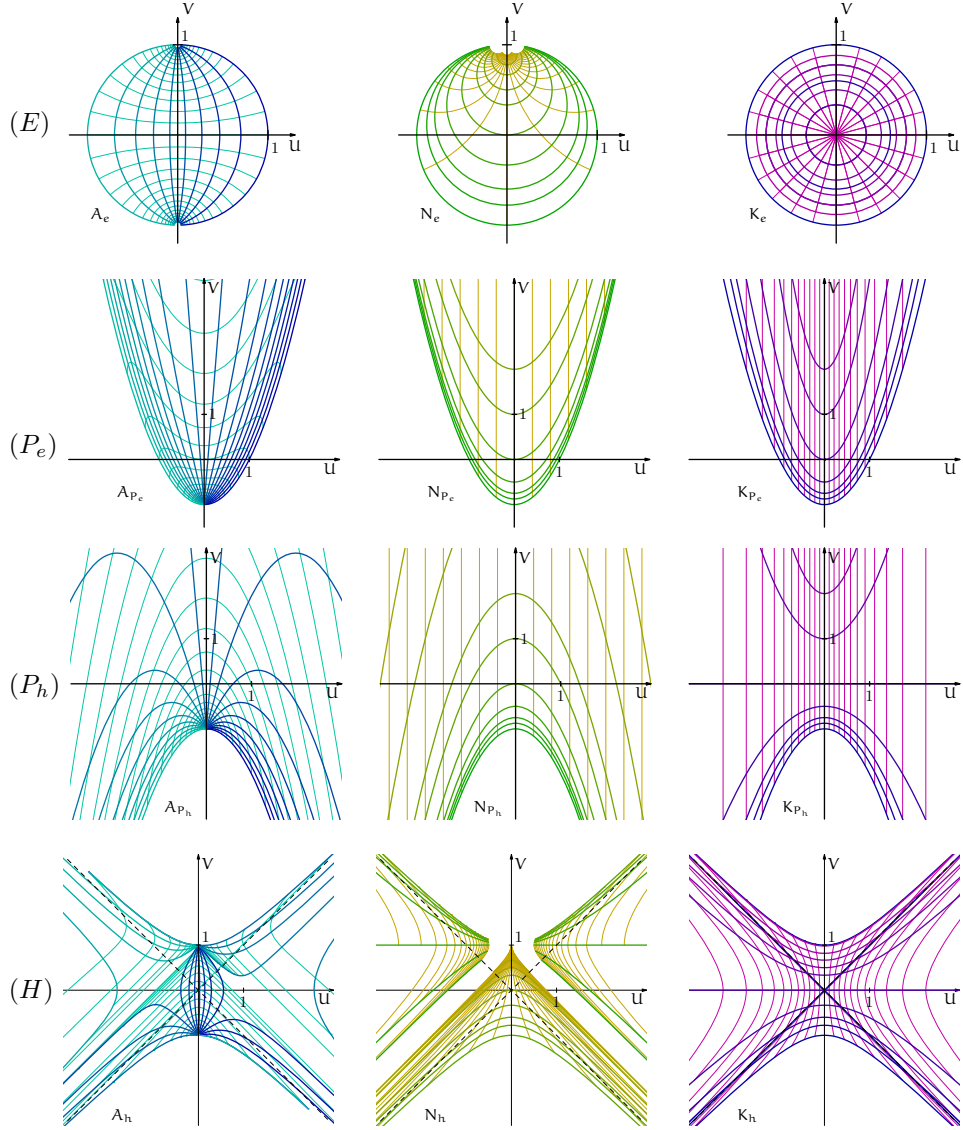

  \centering
\raisebox{1.5cm}{(\(E\))}  \ 
\myincludegraphics[scale=.6
]{\jobname.40}\hfill
  \myincludegraphics[scale=.6
]{\jobname.41}\hfill
  \myincludegraphics[scale=.6
]{\jobname.42}\\[4mm]
\raisebox{1cm}{(\(P_e\))}
\myincludegraphics[scale=.6
]{\jobname.50}\hfill
  \myincludegraphics[scale=.6
]{\jobname.51}\hfill
  \myincludegraphics[scale=.6
]{\jobname.52}\\[2mm]
\raisebox{1.7cm}{(\(P_h\))} 
 \myincludegraphics[scale=.6
]{\jobname.53}\hfill
  \myincludegraphics[scale=.6
]{\jobname.54}\hfill
  \myincludegraphics[scale=.6
]{\jobname.55}\\[4mm]
\raisebox{1.7cm}{(\(H\))}
  \myincludegraphics[scale=.6
]{\jobname.60}\hfill
  \myincludegraphics[scale=.6
]{\jobname.61}\hfill
  \myincludegraphics[scale=.6
]{\jobname.62}\hfill
\caption[The elliptic, parabolic and hyperbolic unit disks]{The
  images of unit disks with orbits of subgroups \(A\),
  \(N\) and \(K\) correspondingly:\\
  (\(E\)): The elliptic unit disk;\\
  (\(P_e\)): The first version of parabolic unit disk with an
  \emph{elliptic type} of Cayley transform (the second---pure
  \emph{parabolic type} (\(P_p\)) transform---is the shift down by
  \(1\) of Figures~\ref{fig:a-n-action}
  and~\ref{fig:k-sungroup}(\(K_p\))). \\
  (\(P_h\)): The third version of parabolic unit disk with a
  \emph{hyperbolic type} of Cayley transform.\\
  \ \qquad(\(H\)): The hyperbolic unit disk. 
}
  \label{fig:unit-disks}
\end{figure}

The images of elliptic actions of subgroups \(A\), \(N\), \(K\) are
given in Figure~\ref{fig:unit-disks}(\(E\)).  The types of orbits can be
easily distinguished by the number of \emph{fixed points on the
  boundary}: two, one and zero correspondingly. In some sense the Cayley
transform swaps complexities: in contrast to on the upper half plane the \(K\)-action
is now simple but \(A\) and \(N\) are not. The simplicity of \(K\)
orbits  is explained by the diagonalisation of matrices:
\begin{equation}
  \label{eq:k-diagonalisation}
   \frac{1}{2}
   \begin{pmatrix}
    1 & - e_2 \\ - e_2 & 1
  \end{pmatrix}
  \begin{pmatrix}
    \cos \phi  & -e_1\sin \phi  \\ - e_1\sin \phi & \cos \phi
  \end{pmatrix}
  \begin{pmatrix}
    1 & e_2 \\  e_2 & 1
  \end{pmatrix} 
  =
  \begin{pmatrix}
    \rme^{\rmi \phi} & 0\\ 0 & \rme^{\rmi \phi}
  \end{pmatrix},
\end{equation}
where \(\rmi =e_1e_2\) behaves as the complex imaginary unit, i.e. \(\rmi^2=-1\).

A hyperbolic version of the Cayley transform was used
in~\cite{Kisil97c}. The above formula~\eqref{eq:cayley-matr} in
\(\Space{R}{h}\) becomes as follows: 
\begin{equation}
  \label{eq:cayley-hyp}
   g_C= 
  \begin{pmatrix}
    f & h \\ h & f
\end{pmatrix}, \textrm{ with }
h =a+d-(b+c)e_2e_1  \textrm{ and }
f = (a-d) e_2+ (c- b)e_1,
\end{equation}
with some subtle differences in comparison
with~\eqref{eq:cayley-elliptic}. 
The corresponding \(A\), \(N\) and \(K\) orbits are given on
Figure~\ref{fig:unit-disks}(\(H\)). However there is an important
difference between the elliptic and hyperbolic cases similar to the one
discussed in subsection~\ref{sec:invar-upper-half}.

\begin{lem}
  \begin{enumerate}
  \item In the \emph{elliptic case} the ``real axis'' \(U\) is transformed
    to the unit circle and the upper half plane---to the unit disk: 
    \begin{eqnarray}
      \label{eq:unit-circle-ell}
      \{(u,v) \such v = 0\} &\to& \{ (u',v') \such l_e(u'e_1+v'e_2)= u'^2+v'^2=1\} \\
      \label{eq:unit-disk-ell}
      \{(u,v) \such v > 0\} &\to& \{(u',v') \such l_e(u'e_2+v'e_2)= u'^2+v'^2<1\}.
    \end{eqnarray}
    On both sets \(\SL\) acts transitively and the unit circle is generated,
    for example, by the point \((1,0)\) and the unit disk is generated by
    \((0,0)\). 
  \item In the \emph{hyperbolic case} the ``real axis'' \(U\) is transformed
    to the hyperbolic unit circle: 
    \begin{equation}
      \label{eq:unit-circle-hyp}
      \{(u,v) \such v = 0\} \to \{ (u',v') \such l_h(u'e_1+v'e_2)= u'^2-v'^2=-1\} 
    \end{equation}
    On the  hyperbolic unit circle \(\SL\) acts transitively and it is
    generated, for example, by point \((0,1)\). 

        \(\SL\) acts also \emph{transitively on the whole complement} \(\{(u',v')
    \such l_e(u'e_2+v'e_2)\neq -1\}\) to the unit circle, i.e. on its
    ``inner'' and ``outer'' parts together.
  \end{enumerate}
\end{lem} The last feature of the hyperbolic Cayley transform can be
treated in a way described in the end of
subsection~\ref{sec:invar-upper-half},  see also
Figure~\ref{fig:hyp-upper-half-plane}(b).  With such an arrangement
the hyperbolic Cayley transform maps 
the ``upper half'' plane from Figure~\ref{fig:hyp-upper-half-plane}(a)
onto the ``inner part'' of the unit disk from Figure~\ref{fig:hyp-upper-half-plane}(b) .

\begin{figure}[htbp]
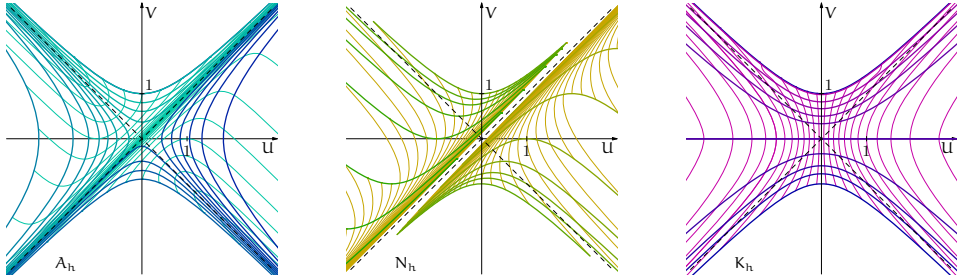

  \centering
  \myincludegraphics[scale=.6
]{\jobname.63}\hfill
  \myincludegraphics[scale=.6
]{\jobname.64}\hfill
  \myincludegraphics[scale=.6
]{\jobname.65}\hfill
  \caption{The hyperbolic unit disk with alternative orbits}
  \label{fig:unit-disks-hyp-alt}
\end{figure}
One may wish that the hyperbolic Cayley transform diagonalises the
action of \(A\) subgroup in a fashion similar
to~\eqref{eq:k-diagonalisation} for \(K\). That is 
achieved~\cite[Ex.~3.1(b)]{Kisil97c} 
if transformation~\eqref{eq:cayley-matr} is preceded by
conjugation with the matrix \(\begin{pmatrix}
  1 & e_1 \\ e_1 & 1
\end{pmatrix}\), or in total the Cayley transform is given by the matrix:
\begin{equation}
  \label{eq:cayle-tr-alt}
  C_1=\begin{pmatrix}
    1-e_2e_1 & e_1-e_2\\
    e_1+\sigma e_2 & 1+ \sigma e_2e_1
  \end{pmatrix} =
  \begin{pmatrix}
    1 & -e_2 \\ \sigma e_2 & 1
  \end{pmatrix}
  \begin{pmatrix}
    1 & e_1 \\ e_1 & 1
  \end{pmatrix}.
\end{equation}
This gives the transformation, cf.~\cite[(3.6--3.7)]{Kisil97c}
\begin{displaymath}
   g_{C_1}= 
  \begin{pmatrix}
    f & h \\ h & f
\end{pmatrix}, \textrm{ with }
h =a(1+e_2e_1)+d(1-e_2e_1)  \textrm{ and }
f = b(e_2-e_1)+c(e_2+e_1),
\end{displaymath}
which obviously keep diagonal form of matrices \(
\begin{pmatrix}
  \alpha^{-1} & 0 \\ 0 & \alpha 
\end{pmatrix}\in A\). 
Orbits of subgroup \(A\), \(N\) and \(K\) for this
transformation are shown on
Figure~\ref{fig:unit-disks-hyp-alt}. The subgroup \(A\) acts by hyperbolic rotations
of \(\Space{R}{h}\) now.

Note that the alternative Cayley
transformation~\eqref{eq:cayle-tr-alt} diagonalises the subgroup \(K\) in
the elliptic case (\(\sigma=-1\)) as well, thus it is also convenient
as an alternative form of the Cayley transform in both the elliptic and
hyperbolic cases. Images of subgroups action in the elliptic case are obtained
from Figure~\ref{fig:unit-disks}(\(E\)) by rotation by \(90^\circ\).

Now we turn to the parabolic case, which benefits from a bigger variety of
choices.  The first natural attempt to define a Cayley
transform can be taken from the same formula~\eqref{eq:cayley-points}
with the parabolic value \(\sigma=0\). The corresponding
transformation defined by the matrix \(
\begin{pmatrix}
  1 & -e_2 \\ 0 & 1
\end{pmatrix}\) turns to be a shift by one unit down. We will denote
it by \(P_p\) for reasons which will become clear shortly.

While being
trivial this transform still possesses some properties of the elliptic
case. For example, the \(K\)-orbits in the elliptic case
(Figure~\ref{fig:unit-disks}(\(K_e\)))  and the
\(A\)-orbits in the hyperbolic case
(Figure~\ref{fig:unit-disks-hyp-alt}(\(A_h\))) are concentric. The
same happens for \(N\)-orbits in the parabolic case
(Figure~\ref{fig:a-n-action}(\(N_a\)))---they all are parabolas
(straight lines) with focus at infinity.

However in the parabolic case it worth considering also both the elliptic \(
\begin{pmatrix}
  1 & e_2 \\ e_2 & 1
\end{pmatrix}\) and hyperbolic \(
\begin{pmatrix}
  1 & e_2 \\ -e_2 & 1
\end{pmatrix}\) transformations~\eqref{eq:cayley-points}. They are
presented on Figure~\ref{fig:unit-disks}, rows (\(P_e\)) and (\(P_h\))
correspondingly. The missing row (\(P_p\)) is formed by the parabolic
transformation discussed in the previous paragraph and illustrated by
Figures~\ref{fig:a-n-action}(\(A_a\)), \ref{fig:a-n-action}(\(N_a\))
and \ref{fig:k-sungroup}(\(K_p\)) with the upper half plane shifted down
by one unit.  Consideration of Figure~\ref{fig:unit-disks} by columns
from top to bottom gives an impressive mixture of many common
properties (e.g. the number of fixed point on the boundary for each subgroup)
with several gradual mutations.

Some properties of parabolic unit disks are as follows:
\begin{lem}
  \label{le:parabolic-disk}
  \begin{enumerate}
  \item All Cayley transforms \(P_e\), \(P_p\) and \(P_h\) act on the
    axises  \(V\) as the shift down by \(1\).
  \item \label{it:parab-unit-disk}
    The parabolic unit disk at \(P_e\) is given by the inequality 
    \(l_p(ue_1+ve_2)\leq 1\) with boundary given by the parabolic unit
    circle \(l_p(ue_1+ve_2)= 1\) in sharp resemblance
    to~\eqref{eq:unit-circle-ell} and~\eqref{eq:unit-circle-hyp}.
  \item The parabolic unit disk at \(P_h\) is given by the inequality 
    \(l_p(-(u+2)e_1-ve_2)\leq 1\) with boundary given by the parabolic unit
    circle \(l_p(-(u+2)e_1+ve_2)= 1\).
  \item\label{it:parabolic-disk3}
    \(N\)-orbits in both transforms \(P_e\) and \(P_h\)
    are \emph{parabolas} with focal length \(1/2\).
  \item \(A\)-orbits in transforms \(P_e\) and \(P_h\) are segments of
    parabolas with focal length \(1/2\) 
    passing through \((0,-1)\). Their vertices belongs to two parabolas
    \(v=-x^2-1\) and \(v=x^2-1\) correspondingly, which are boundaries
    of parabolic circles in  \(P_h\) and \(P_e\) (note the swap!).
  \item\label{it:parabolic-disk5}
    \(K\)-orbits in transform \(P_e\)  are parabolas
    with focal length less than \(1/2\) and in transform
    \(P_h\)---with inverse of focal length bigger than \(-2\).
  \end{enumerate}
\end{lem}
Of course property~\ref{it:parab-unit-disk} makes transformations
\(P_e\) very appealing as a ``right'' parabolic version of the Cayley
transform. However it seems that all three transformations
\(P_{e,p,h}\) have their own merits which may be decisive in
particular circumstances.

\begin{rem}
  \label{re:par-more-cayley}
  We see that the varieties of possible Cayley transforms in the parabolic
  case is bigger than in the two other cases. It is interesting that this parabolic
  richness is a consequence of the parabolic degeneracy of the generator
  \(e_2^2=0\). Indeed for both the elliptic and the hyperbolic sign in \(e_2^2=\pm 1\) only
  one matrix~\eqref{eq:cayley-points} out of two possible \(
  \begin{pmatrix}
    1 & -e_2 \\ \pm e_2 & 1
  \end{pmatrix}\)\vspace{3pt} has a non-zero determinant. And only for
  the degenerate parabolic value \(e_2^2=0\) both these matrices are
  non-degenerate!

\end{rem}

\section*{Acknowledgments}
\label{sec:acknowledgments}
The first named author is grateful to Prof.~S.~Plaksa for pointing out
the relevance of the book~\cite{LavrentShabat77} to the present paper.
The first named author thanks Prof.~S.Blumin for
informing us on the highly relevant book~\cite{Yaglom79} after this
paper was distributed as a preprint. Last but not least we are grateful
to Dr~I.R.~Porteous for the careful reading of the paper and
numerous comments and remarks.

The extensive graphics in this paper were produced with the help of
\GiNaC~\cite{GiNaC} computer algebra system. Since this tool is of separate
interest we explain its usage by examples from this article in the
separate paper~\cite{Kisil04c}. The \NoWEB~\cite{NoWEB}  wrapper for \CPP\ source
code is included in the \texttt{arXiv.org} files of both these
papers~\cite{Kisil04c,Kisil04b}.

\small
\bibliographystyle{plain}
\bibliography{arare,abbrevmr,akisil,analyse,aphysics}
\end{document}